\magnification\magstep1

\magnification\magstep1

\def\sqr#1#2{{\vcenter{\vbox{\hrule height.#2pt
    \hbox{\vrule width.#2pt height#1pt \kern#1pt
    \vrule width.#2pt}
    \hrule height.#2pt}}}}

\def\ecarre{\; \mathchoice\sqr56\sqr56\sqr{4.3}5\sqr{3.7}5}

\def\Ga{\Gamma}
\def\ga{\gamma}
\def\al{\alpha}

\def\de{\delta}

\def\om{\omega}

\def\pa{\partial}

\def\tr{\hbox{tr}}

\font\bbb=msbm10

\def\R{\hbox{{\bbb R}}}

\def\H{\hbox{{\bbb H}}}

\def\ssni{\smallskip\noindent}
\def\msni{\medskip\noindent}

\def\ni{\noindent}
\def\ss{\smallskip}
\def\ms{\medskip}
\def\bs{\bigskip}

\centerline{\bf FEDOSOV MANIFOLDS}

\bs
\centerline{\bf Israel Gelfand}
\centerline{(Rutgers University)}
\centerline{igelfand@math.rutgers.edu}
\ms
\centerline{\bf  Vladimir Retakh}
\centerline{(University of Arkansas, Fayetteville)}
\centerline{retakh@math.harvard.edu}
\ms
\centerline{\bf Mikhail Shubin}
\centerline{(Northeastern University)}
\centerline{shubin@neu.edu}
\bigskip
\centerline{\it Presented at Symplectic Geometry Workshop, Toronto,
June 1997}
\bigskip
\centerline{\bf 1. Introduction.}

\medskip
{\bf A.} Connection (or parallel transport) is the main object in the
geometry of manifolds.
Indeed, connections allow comparison between geometric quantities
associated with different
(distant) points of the manifold. In 1917 Levi-Civita [Lv] introduced
a notion of connection for a manifold embedded into $\R^n$. In
1918 H. Weyl [We] introduced general symmetric linear connections
in the tangent bundle. In
1922 E. Cartan [C1], [C2], [C3] studied non-symmetric linear connections
(or affine connections),
 where the torsion is interpreted in general relativity as the density
of intrinsic angular momentum (spin) of matter. Actually Cartan explains
that affine connection and linear connection with torsion is the same
thing.
According to Cartan, connection
is a mathematical alias for an observer traveling in space-time and carrying
measuring instruments. Global comparison devices (like a transformation
group in the  Erlangen program by F.Klein) are forbidden by the finite
propagation
speed restriction (see, for example, [S]).

\ms
It might seem that the Riemannian geometry only requires a Riemannian metric,
but this is due to the fact that there is a canonical Levi-Civita connection
associated with this metric.

\ms
In this paper we  study geometry of {\it symplectic connections},
i.e. symmetric (torsion-free) connections which preserve a given symplectic
form.

\ms
We start (in Sect.1) with a more general context. Let us consider an
{\it almost symplectic manifold}, i.e. manifold $M$ with a non-degenerate
exterior 2-form $\om$ (not necessarily closed) and try to describe all linear
connections (i.e. connections in the tangent bundle) which  preserve this
form.
In the classical  case of non-degenerate symmetric form $g=(g_{ij})$
(Riemannian or pseudo-Riemannian) the answer
is well known due to Levi-Civita: for any  tensor $T=(T^k_{ij})$, which is
antisymmetric with respect to the subscripts $i,j$, there exists a unique
connection
$\Ga=(\Ga^k_{ij})$ with the torsion tensor $T$ (i.e.
$\Ga^k_{ij}-\Ga^k_{ji}=T^k_{ij}$ )
such that $\Ga$ preserves $g$. If $T=0$, this connection is called
the {\it Levi-Civita connection}. It is the {\it canonical} connection
associated
with $g$.

\ms
We first prove a skew-symmetric analog of this Levi-Civita theorem. Namely,
there is a canonical one-one correspondence between the following two sets
of linear connections on an almost symplectic manifold $(M,\om)$:

\ms
1) all symmetric connections on $M$;

\ss
2) all connections on $M$ which preserve $\om$.

\msni
This correspondence is obtained as follows:
for a connection $\Ga$ which preserves $\om$, we  just take its symmetric part
$\Pi=(\Pi^i_{jk})$, $\Pi^i_{jk}={1\over 2}(\Ga^i_{jk}+\Ga^i_{kj})$
(which is again a connection).
This gives a bijection between all  connections preserving $\om$ and all
symmetric connections. The inverse operator is also given by explicit formulas
(see (1.7) below).

\medskip
Note the following analogy with the  Levi-Civita
construction. There we can prescribe the antisymmetric part of the
connection, which is a tensor (the torsion tensor). In particular we can take
this tensor to be zero (then we arrive to the Levi-Civita connection).
In symplectic case (or rather in case of any non-degenerate antisymmetric
form) we can
prescribe the symmetric part of the connection which is not a tensor
but a symmetric connection.

\ms
Geometry of almost symplectic manifolds was studied by H.-C.Lee [Le1], [Le2],
P.Li\-ber\-mann [Li1-Li3], V.G.Lemlein [Lm],
Ph.Tondeur [To] and I.Vaisman [Va1] (see also references in [Li2, Li3], [Va1]).
In particular, V.G.Lemlein and Ph.Tondeur independently gave two different
(but less geometrical) description of
all connections preserving $\om$. (See also [Va1] for a proof and a
refinement of the
Ph.Tondeur result.)

\ms
V.G.Lemlein [Lm] also introduced almost symplectic manifolds with a
symmetric connection
such that $\nabla\om=\mu(x)d\om$ (as tensors). It follows that
$\mu(x)=1/3$. He described
the connections  $\nabla$ with this property in terms of an object
$\ga_{ijk}$ which is
completely symmetric in $i,j,k$ (though it is neither a tensor nor a
connection).

\ms
{\bf B.} In Sect.2 we introduce our main object: {\it Fedosov manifolds}.
By definition a Fedosov manifold is a triple $(M,\om,\Ga)$ where $(M,\om)$ is
a symplectic manifold (i.e.  $\om$
is a symplectic form, a non-degenerate closed exterior 2-form, on a
$C^\infty$-manifold $M$)
and $\Ga$ is a
symplectic connection on $M$.

\ms
The famous result of Fedosov gives a canonical deformation quantization on
such a manifold.
This quantization is canonically defined by the data $(M,\om,\Ga)$ (see
e.g.
[F1], [F2], [DW-L], [W], [D]).

\ms
For the role of symplectic connections in geometric quantization
see e.g. [Lch], [Hs]. A symplectic connection can be uniquely associated to
a symplectic
manifold equipped with two transversal polarizations (involutive Lagrangian
distributions
in the complexified tangent bundles) --  see [Hs], [M], [Va1]. K.Habermann
[Ha] defined and
studied symplectic Dirac operators on Fedosov manifolds with an additional
metaplectic
structure.

\ms
Note that symplectic connections exist on  any symplectic manifold. Indeed,
locally such a connection can be constructed in Darboux coordinates by
assigning
all Christoffel symbols to be identically $0$. Then we can glue a global
symplectic connection by a partition of unity.  Therefore
a structure of a symplectic manifold can be always extended to a structure
of a Fedosov manifold.

\ms
In fact the symplectic connections on a symplectic manifold form an
infinite-di\-men\-sio\-nal
affine space with difference vector-space isomorphic to the space of all
completely
symmetric 3-covariant tensors - see e.g. [BFFLS], [Ve] and Section 1
(the paragraph after (1.5)) below.

\ms
Fedosov manifolds constitute a natural generalization of K\"ahler manifolds.
A K\"ahler manifold can be defined as a Riemannian  almost complex manifold
$(M,g,J)$
(here $g$ is a Riemannian metric and $J$ an almost complex structure
on the manifold $M$), such that  $J$ is orthogonal with respect to $g$  and
parallel with respect to the Levi-Civita connection $\Ga$ associated with $g$
(see e.g. [K-N], vol.2, Ch.IX, Sect.4, or [Br], p.148).
The last property can be formulated
by saying that the 2-form $\om$ defined by $g,J$ according to the formula
$$\om(X,Y)=g(X,JY)\;,\leqno (0.1)$$
is $\Ga$-parallel, i.e. invariant under the parallel transport defined by
$\Ga$ (or,
equivalently, its covariant derivative vanishes).
It follows that $\om$ is closed, hence symplectic because it is
automatically non-degenerate.
Therefore every K\"ahler manifold is a Fedosov manifold.

\ms
On the other hand there are symplectic manifolds such that they do not admit
any  K\"ahler structure, even with a different $\om$ (see e.g. [Br], p.147,
for a Thurston example).

\ms
We study the curvature tensor
of a symplectic connection (on a Fedosov manifold). The usual curvature
tensor $R^i_{\ jkl}$
is well defined by the connection alone.
But the tensor $R_{ijkl}=\om_{ip}R^p_{\ jkl}$ is more interesting. (We use
the standard Einstein summation convention, so here the summation over $p$
is understood.)
This tensor is well known. For example, V.G.Lemlein [Lm] introduced it in a
more general context
of almost symplectic manifolds with a connection, I.Vaisman [Va1] obtained
important
symmetry properties. Since we tried to make this paper self-contained, we give
(usually new) proofs
of some properties of this tensor which were first obtained in [Va1]. Following
I.Vaisman we will call $R_{ijkl}$ {\it the  symplectic curvature tensor}.

\ms
The symplectic curvature tensor has  symmetries which are different from
the symmetries of the
curvature tensor on a Riemannian manifold. Namely, aside of standard symmetries
$R_{ijkl}=-R_{ijlk}$, $R_{i(ijk)}=0$
(here the parentheses
mean that we should take sum over all cyclic permutations of the subscripts
inside the parentheses),
it is symmetric with respect to the first two subscripts: $R_{ijkl}=R_{jikl}$.
We deduce from these identities that $R_{(ijkl)}=0$.

\ms
{\bf C.} In Sect.3 we consider the Ricci tensor $K_{ij}=R^k_{\
ikj}=\om^{kp}R_{pikj}$
which is again defined
by the connection alone. I.Vaisman [Va1] proved that it is in fact symmetric
(as in the case of Riemannian manifolds)
and the structure of the Fedosov manifold plays important role in the proof.
We give two proofs of the  symmetry of $K_{ij}$ (in particular we deduce it
from the identity $R_{(ijkl)}=0$).

\ms
There are two natural non-trivial contractions of the curvature tensor
on a Fedosov manifold seemingly leading to two different Ricci type tensors
(tensors of type (2,0) or bilinear forms on each tangent space).
However we prove that they only differ by a constant factor 2:
$\om^{kl}R_{ijkl}=2K_{ij}$. We also  prove that the operator
$$L=\om^{ij}\nabla_i\nabla_j\;,$$
on vector fields has order zero and corresponds to the Ricci tensor
(with one subscript lifted by $\om$).

\ms
In the Riemannian case contracting the Ricci tensor we obtain the scalar
curvature.
On a Fedosov manifold we can only use the symplectic form for the
contraction, and this
obviously leads to a trivial result.

\ms
Note also that I.Vaisman [Va1] provided a complete group-theoretic analysis
of the space
of all curvature-type tensors  (i.e. tensors with the symmetries described
above) as
a representation space of the group $\hbox{Sp}(2n)$, where $2n=\dim_{\R}
M$. He proved that
the corresponding representation space splits into a direct sum of two
suspaces,
one of them is the space of all tensors with vanishing Ricci curvature
part, and another
is also explicitly described. Each of these two subspaces is irreducible.

\ms
Vanishing of the Ricci tensor  provides  the vacuum Einstein equations in the
Riemannian or pseudo-Riemannian case. This is a determined system for the components 
of the metric.
We can also write the vacuum Einstein equations $K_{ij}=0$
on a Fedosov manifold.
We call the manifolds with this property {\it symplectic Einstein manifolds.}
(They are called Ricci flat in [Va1].)
However $K_{ij}=0$ is an undetermined system
(for the symplectic form together with the symplectic connection). The obvious
local gauge transformations (symplectic transformations of coordinates) are not
sufficient to remove the underdeterminacy even in dimension 2.

\ms
We also give a definition of sectional curvature on a Fedosov manifold. It
is a function
on non-isotropic 2-planes $\Pi$ in the tangent bundle, and it is defined  by
the canonical form of the quadratic form $Z\mapsto R(Z,Z,X,Y)$, $Z\in\Pi$,
where $\{X,Y\}$ is a basis of $\Pi$ such that $\om(X,Y)=1$. The quadratic
form may be
degenerate of rank $0$ or $1$. If it is non-degenerate, then it has a real
(non-zero)
invariant and can also be elliptic or hyperbolic.

\ms
{\bf D.} In Sect.4 we introduce {\it normal coordinates} on a Fedosov
manifold. They are defined
as affine coordinates in the tangent space to $T_pM$ which are transferred
to $M$
by the exponential map of the given connection $\Ga$ at $p$.
Since they are defined by the connection only, the construction does not
differ from the
Riemannian case. We can also define affine {\it extensions} of an arbitrary
tensor
$T=(T^{i_1\dots i_k}_{j_1\dots j_l})$ which are tensors
$T^{i_1\dots i_k}_{j_1\dots j_l\al_1\dots\al_r}$  obtained at $p\in M$ by
taking
the corresponding derivatives of the components of $T$ in normal
coordinates centered at $p$
and then evaluating the result at $p$. Applying this procedure to the
Christoffel symbols $\Ga_{ijk}$, we obtain so called {\it normal tensors}
$$A_{ijk}=0,\  A_{ijk\al_1},\  A_{ijk\al_1\al_2},\dots\;. $$
(They are tensors in spite of the fact that $\Ga_{ijk}$ do not form a tensor,
and due to the fact that $\Ga_{ijk}$ undergo the tensor transformation law
if the
transformation of the coordinates is linear.)

Extending  results of T.Y.Thomas and O.Veblen (see e.g. [E, Th, V]) we
describe all
algebraic relations  for the normal tensors and the extensions of the
symplectic form at a point of a Fedosov manifold. Compared with the classical
situation (arbitrary symmetric connection) a new relation is added to provide
the compatibility of the normal tensors with the existence of the preserved
symplectic form (see (4.18)).
 We also prove that all local
invariants of a Fedosov manifold are appropriate functions of the
components of
the normal tensors and of the components of the symplectic form.

\ms
{\bf E.} In Section 5 we study relations between the normal tensors and the
curvature tensor.
Here we again use an appropriately modified technique by T.Y.Thomas and
O.Veblen.
We provide formulas which express  normal tensors through the extensions of
the curvature tensor (or its covariant derivatives) and vice versa. A
simplest corollary
of these formulas is
$$\om_{ij,kl}={1\over 3}R_{klij}\;,$$
where $\om_{ij,kl}$ is the 2nd extension of $\om$ and $R$ is the curvature
tensor.

\ms
We also prove that the above mentioned symmetries of the curvature tensor
together
with the Bianchi identity (see (5.25)) and an integrability identity (see
(5.26))
provide a complete system of identities
for the curvature tensor and its first covariant derivatives at a point of
a Fedosov
manifold.

\ms It follows from the results of Sect.5 that
 any local invariant of a Fedosov manifold is an appropriate function
of the components of $\om$ and of the covariant derivatives of the curvature
tensor (Theorem 5.11).
A complete description of such functions in Darboux coordinates  was done
by D.E.Tamarkin [Ta] . We plan to discuss local
invariants in more detail in a future paper.

\ms
{\bf F.} The Poisson manifolds are a natural generalization of symplectic
manifolds
(see e.g. [Va2]).
They only have $\om^{ij}$ which can be degenerate, but satisfy a natural
identity
which is equivalent to the Jacobi identity for  the corresponding Poisson
bracket.
Recently M.Kontsevich [Ko] constructed a deformation quantization on any
such manifold.
Note however that a connection preserving the Poisson structure does not
always exist. Indeed,
if it exists then the matrices $\om^{ij}$ should  have a locally constant
rank,
because the
parallel transport  preserves the rank. This is not the case
e.g. if all the components $\om^{ij}(x)$ are linear in $x$ and not
identically $0$,
as is the case for the linear Poisson brackets defined by non-abelian Lie
algebras.

\ms
{\bf G.} We are grateful to M.Eastwood, S.Gelfand, A.R.Gover,  M.Khalkhali,
J.Slovak, P.Tod
and I.Vaisman for useful discussions and references. The third author also
acknowledges support
of the University of Adelaide (Australia) and the Fields Institute (Toronto)
during the final stage  of his work on this paper.

\bigskip
\centerline{\bf 1. Connections on almost symplectic manifolds.}
\medskip
Suppose we are given a manifold $M,\ \dim M=n,$ and a non-degenerate 2-form
$\omega$
which we shall fix. Such a pair $(M,\om)$ is called {\it almost symplectic
manifold}
(see e.g. [Li] and references there).

Let $\nabla$ be  a connection (covariant derivative) on $M$. We will say that
$\nabla$ {\it preserves} $\omega$ if $\nabla\omega=0$ or
$$Z(\omega(X,Y))=\omega(\nabla_Z X,Y)+\omega(X,\nabla_Z Y) \leqno(1.1)$$
for any vector fields $X,Y,Z$.

\ms
Let us write this in coordinates. If $x^1,\dots,x^n$ are local coordinates,
denote $\partial_i={\partial\over\partial x^i}$. The components of $\omega$ in
these coordinates are $\omega_{ij}=\omega(\partial_i, \partial_j)$. Let us
denote $\nabla_i=\nabla_{\partial_{i}}$ and introduce the {\it Christoffel
symbols}
$\Gamma^k_{ij}$ by $\nabla_i\partial_j=\Gamma^k_{ij}\partial_k$.
(Here and later on we  omit the sign of summation using the standard
Einstein convention.) We will identify the connection $\nabla$ with the set
$\Ga=(\Ga^k_{ij})$
of its Christoffel symbols.

\ms
It is
sufficient to write (1.1) for $X=\partial_i,\  Y=\partial_j,\ Z=\partial_k$.
This gives
$$\pa_k\omega_{ij}=\omega_{il}\Gamma^l_{kj}-\omega_{jl}\Gamma^l_{ki}=
\Gamma_{ikj}-\Gamma_{jki}\leqno(1.2)$$
where $\Gamma_{ikj}=\omega_{il}\Gamma^l_{kj}$.

\ms
The equality $d\omega=0$ or $d(\omega_{ij}dx^i\wedge dx^j)=0$ means
$$\pa_k\omega_{ij}+\pa_j\omega_{ki}+\pa_i\omega_{jk}=0\,.\leqno(1.3)$$

Recall that $\nabla$ is symmetric
(or has no torsion) if
$$\nabla_X Y-\nabla_Y X=[X,Y]\leqno(1.4)$$
for arbitrary vector fields $X,Y$.
Taking $X=\partial_i$, $Y=\partial_j$, we can rewrite this
in coordinates  as $\Gamma_{ij}^k=\Gamma^k_{ji}$ or, equivalently,
$$\Gamma_{ijk}=\Gamma_{ikj}\,.\leqno(1.5)$$

\ms
If $d\om=0$ then we can locally find
Darboux coordinates where $\omega_{ij}$ are constants. Then (1.2) is
equivalent to
the symmetry of $\Gamma_{ijk}$ with respect to $i,k$. Therefore in Darboux
coordinates the symmetric connections preserving $\omega$ are exactly the
connections with the Christoffel symbols $\Gamma_{ijk}$ which are completely
symmetric with respect to all indices $i,j,k$. Now note that the difference
between two connections is a tensor. It follows that for any two symmetric
connections
$\nabla,\ \nabla^\prime$
preserving the same form $\omega$,  the difference
$\Gamma_{ijk}-\Gamma^\prime_{ijk}$ is
a 3-covariant tensor which is completely symmetric in all indices in
Darboux coordinates, hence in any coordinates.

\ms
Let $\Gamma=(\Gamma_{ij}^k)$ be a linear connection on $M$. Denote by
$\Pi=(\Pi_{ij}^k)$ its
symmetric part ($\Pi=\Pi(\Gamma)$) i.e.
$$\Pi_{ij}^k={1\over 2}(\Gamma_{ij}^k+\Gamma_{ji}^k).\leqno(1.6)$$
\medskip
{\bf Theorem 1.1.} {\it The map} $\Gamma\mapsto\Pi(\Gamma)$ {\it gives a
bijective
affine correspondence between the set of all connections preserving}
$\omega$ {\it and
the set of all symmetric connections. The inverse map} $\Pi\mapsto\Gamma$
{\it is
given by}
$$\Gamma_{kij}={1\over 2}(\pa_k\omega_{ij}-\pa_i\omega_{jk}-\pa_j\omega_{ki})+
(\Pi_{kij}+\Pi_{jik}-\Pi_{ijk}).\leqno(1.7)$$
{\it In case when} $\omega$ {\it is closed (hence symplectic) this formula
can be
rewritten as}
$$\Gamma_{kij}=\pa_k\omega_{ij}+(\Pi_{kij}+\Pi_{jik}-\Pi_{ijk}). \leqno(1.8)$$

\medskip
{\bf Proof.}  Assume that $\Ga$ preserves  $\omega$ i.e. (1.2) is satisfied.
It follows by the cyclic permutation of $i,j,k$ that
$$\pa_i\omega_{jk}=\Gamma_{jik}-\Gamma_{kij}$$
and
$$\pa_j\omega_{ki}=\Gamma_{kji}-\Gamma_{ijk}.$$
Subtracting these  two relations from (1.2) and using (1.6) we come to
(1.7). Now (1.8) follows if we use (1.3). $\ecarre$

\ms
{\bf Corollary 1.2.}
{\it If} $\Gamma$ {\it is a connection which preserves} $\omega$ {\it then the
following properties are equivalent:}

(i) $\Gamma$ {\it is symmetric (or $\Ga=\Pi$)};

(ii) $\Pi$ {\it preserves $\omega$, i.e. the symmetric part of $\Ga$
preserves $\om$.}
\ms
{\bf Proof.}
Clearly (i) implies (ii).
Vice versa, assume that (ii) is true. Then (i) should be true because a
connection which
preserves $\omega$ and has the given symmetric part $\Pi$ is unique due to
the theorem.
$\ecarre$

\ms
{\bf Corollary 1.3.} [Va1] {\it Let $(M,\om)$ be an almost symplectic
manifold with an action
of a group $G$ which preserves $\om$. Assume that there
exists a $G$-invariant connection on $M$. Then $M$ has a $G$-invariant
connection
preserving $\om$.}

\ms
{\bf Proof.}  If $\Ga$ is a $G$-invariant connection on $M$, then it
symmetric part $\Pi$
is also $G$-invariant. Therefore the corresponding $\om$-preserving
connection (from Theorem 1.1)
will be also $G$-invariant. $\ecarre$

\bigskip
{\bf Remark 1.4.} It is a well known fact that a symmetric connection
preserving
$\omega$ exists if and only
if $\omega$ is closed (see e.g. [Li1], [To], [Va1]). This can be proved as
follows: if $\Ga$ is
a connection, which preserves $\om$, i.e.
(1.2) is satisfied, then we can rewrite (1.3) in the form
$$\Gamma_{ikj}-\Gamma_{jki}+\Gamma_{kji}-\Gamma_{ijk}+
\Gamma_{jik}-\Gamma_{kij}=0$$
which is obviously true for any symmetric connection.

\ms
Vice versa, if $\om$ is closed (hence symplectic), then locally we can take
the trivial connection in Darboux coordinates. Globally we can glue a symmetric
connection preserving $\om$ using a partition of unity.

\medskip
{\bf Remark 1.5.} Suppose we are given non-degenerate symmetric and
antisymmetric
2-forms $g$ and $\omega$ at the same time. Then we can construct a sequence of
connections. Let us
start with any antisymmetric torsion tensor $T_0$ (e.g. zero). Then there
exists
a unique $g$-preserving  connection $\Pi_1^g$
with the torsion $T_0$. Denote by $S_1$ its symmetric part (this is a symmetric
connection). Then we can find a (unique) $\omega$-preserving connection
$\Gamma_1^\omega$ which has the symmetric part $S_1$. Denote its torsion tensor
by $T_1$. Now we can repeat the construction starting with $T_1$ instead of
$T_0$.
In this way we obtain a sequence of connections and torsion tensors
$$T_0\mapsto\Pi_1^g\mapsto S_1\mapsto\Gamma_1^\omega\mapsto T_1\mapsto
\Pi_2^g\mapsto S_2\mapsto\Gamma_2^\omega\mapsto T_2\mapsto\dots\;.$$
This sequence looks similar to the sequence which appears when one uses two
symplectic
structures in the theory of integrable systems.

Note that for K\"ahler manifolds, starting with $T_0=0$ we obtain the sequence
which stabilizes i.e. we get $T_j=0$, $\Pi_j^g=S_j=\Gamma_j^\omega$ for all
$j$.

\ms
{\bf Remark 1.6.}
Independently of  Theorem 1.1 it is easy to observe a coincidence of
functional dimensions (the number of independent functional parameters) of the
two sets of connections above: all connections preserving $\om$ and all
symmetric  connections.  Playing with the dimensions we noticed
another curious coincidence:
$$C_n-C_n^\omega=S_n-S_n^\omega+\Lambda_n^3, \leqno (1.9)$$
where $\omega$ is a given non-degenerate (exterior) 2-form,
$$C_n=\hbox{f-dim of all connections}=n^3,$$
(f-dim means functional dimension)
$$S_n^\omega=\hbox{f-dim of symmetric connections preserving\ }\omega=
{(n+2)(n+1)n\over 6},$$
$$C_n^\omega=\hbox{f-dim of connections preserving\ }\omega={n^2(n+1)\over
2},$$
$$S_n=\hbox{f-dim of all symmetric connections}={n^2(n+1)\over 2},$$
$$\Lambda_n^3=\hbox{f-dim of all (external) 3-forms}={n(n-1)(n-2)\over 6}.$$
It would be interesting to clarify the coincidence (1.9)
(similarly to the  clarification of the
``coincidence'' $C_n^\omega=S_n$ which was given above).

\bs
\centerline{\bf 2. Symplectic connections and their curvature.}
\ms
Let $(M,\om)$ be a symplectic manifold, $\nabla$ (or $\Ga$) a connection on
$M$.
\ms
{\bf Definition 2.1.} If $\nabla$ is  symmetric and preserves
the given symplectic form $\omega$ then we will say that $\nabla$ is a
{\it symplectic connection}.
\medskip
{\bf Definition 2.2.} {\it Fedosov manifold} is a symplectic manifold with a
given symplectic connection.
A Fedosov manifold $(M,\om,\Ga)$ is called {\it real-analytic} if
$M,\om,\Ga$ are
all real-analytic.
\ms
The curvature tensor of a symplectic connection is defined
by the usual formula
$$R(X,Y)Z=\nabla_X\nabla_Y Z-\nabla_Y\nabla_X Z-\nabla_{[X,Y]}Z.\leqno(2.1)$$
The components of the curvature tensor are introduced by
$$R(\partial_j,\partial_k)\partial_i=R^m_{\ ijk}\partial_m\,.\leqno(2.2)$$
Denote also
$$R_{ijkl}=\omega_{im}R^m_{\ jkl}=\omega(\partial_i,
R(\partial_k,\partial_l)\partial_j)\,.\leqno(2.3)$$
Instead of $R_{ijkl}$ we can also consider $R(X,Y,Z,W)$  which is a
multilinear function
on any tangent space $T_xM$ (or a function of 4 vector fields which is
multilinear over
$C ^\infty(M)$):
$$R(X,Y,Z,W)=\om(X,R(Z,W)Y), \leqno (2.3')$$
so that $R_{ijkl}=R(\pa_i,\pa_j,\pa_k,\pa_l)$. But it is usually more
convenient to do
calculations using components.

\msni
It is obvious that
$$R_{ijkl}=-R_{ijlk}\,.\leqno(2.4)$$
The formula expressing the components of the curvature tensor in terms of the
Christoffel symbols has the standard form:
$$R^l_{\ ijk}=\partial_j\Gamma^l_{ki}-\partial_k\Gamma^l_{ij}+
\Gamma^m_{ki}\Gamma^l_{jm}-\Gamma^m_{ij}\Gamma^l_{km}\,.\leqno(2.5)$$
For any symmetric connection we have the first Bianchi identity
$$R_{i(jkl)}:=R_{ijkl}+R_{iklj}+R_{iljk}=0\,.\leqno(2.6)$$
( For the proofs of (2.5) and (2.6) see e.g.
[D-F-N], Sect.30, or [He], Ch.1, Sect.8, 12).
\bigskip
{\bf Proposition 2.3.} [Va1] {\it For any symplectic connection}
$$R_{ijkl}=R_{jikl}\,.\leqno(2.7)$$
\medskip
{\bf Proof.} Let us consider
$$\pa_k\pa_l\omega_{ij}=
\pa_k[\omega(\nabla_l\pa_i,\pa_j)+
\omega(\pa_i,\nabla_l\pa_j)]=\leqno(2.8)$$
$$=\omega(\nabla_k\nabla_l\pa_i,\pa_j)+
\omega(\nabla_l\pa_i,\nabla_k\pa_j)+
\omega(\nabla_k\pa_i,\nabla_l\pa_j)+
\omega(\pa_i,\nabla_k\nabla_l\pa_j)\,.$$
Changing places of $k,l$ and subtracting the result from the (2.8) we obtain:
$$0=\omega([\nabla_k,\nabla_l]\pa_i,\pa_j)+
\omega(\pa_i,[\nabla_k,\nabla_l]\pa_j)=
\omega(R^m_{\ ikl}\pa_m,\partial_j)+\omega(\pa_i,R^m_{\ jkl}\pa_m)=$$
$$=\omega_{mj}R^m_{\ ikl}+\omega_{im}R^m_{\ jkl}=-R_{jikl}+R_{ijkl}\,.\ecarre$$
\bigskip
{\bf Proposition 2.4.} {\it For any symplectic connection}
 $$R_{(ijkl)}:=R_{ijkl}+R_{lijk}+R_{klij}+R_{jkli}=0\,.\leqno(2.9)$$
\medskip
{\bf Proof.} Using (2.6) we obtain
$$R_{ijkl}+R_{iklj}+R_{iljk}=0\,,$$
$$R_{jikl}+R_{jkli}+R_{jlik}=0\,,$$
$$R_{kijl}+R_{kjli}+R_{klij}=0\,,$$
$$R_{lijk}+R_{ljki}+R_{lkij}=0\,.$$
Adding all these equalities and using (2.4) and (2.7), we obtain (2.9).
$\ecarre$

\bs
\centerline{\bf 3. Ricci tensor and sectional curvature.}
\ms
Let us look for possible contractions of the curvature tensor which can
lead to new
non-trivial tensors.
\medskip
{\bf Lemma 3.1.} {\it On any Fedosov manifold}
$$R^k_{\ kij}=\omega^{kl}R_{lkij}=0\,.\leqno(3.1)$$
\medskip
{\bf Proof.} The result immediately follows from the symmetry of $R_{lkij}$
with
respect to the first two subscripts (Proposition 2.3). $\ecarre$
\bigskip
{\bf Definition 3.2.}  The {\it Ricci tensor} of a Fedosov manifold is defined
by the formula
$$K_{ij}=R^k_{\ ikj}=\omega^{kl}R_{likj}\,.\leqno(3.2)$$
Equivalently we can define $K$ as a bilinear form on $T_xM$:
$$K(X,Y)=K_{ij}X^iY^j=\tr\{Z\mapsto R(Z,Y)X\}\;,\leqno(3.2')$$
where $X,Y,Z\in T_xM$, $X=X^i\pa_i, Y=Y^j\pa_j$.

\ms
Note that the Ricci tensor depends on connection only. (It does not depend
on $\om$.)
\ms
{\bf Proposition 3.3.} [Va1] {\it On any Fedosov manifold the Ricci tensor
is symmetric:}
$$K_{ij}=K_{ji}\,.\leqno(3.3)$$
\medskip
{\bf Proof.} Multiplying (2.9) by $\omega^{ki}$, summing over $k,i$ and using
the symmetries (2.4) and (2.7) of the curvature tensor we get the desired
result.
$\ecarre$
\bigskip
{\bf Corollary 3.4.} [Va1] {\it On any  Fedosov manifold}
$$\omega^{ji}K_{ij}=0\,.\leqno(3.4)$$

This means that the usual way to define the scalar curvature gives a trivial
result.

\ms
Now let us show that the second non-trivial contraction of the curvature tensor
again leads to the Ricci tensor.

\ms
{\bf Proposition 3.5.} {\it On any Fedosov manifold}
$$\om^{kl}R_{ijkl}=2K_{ij}\;.\leqno (3.5)$$
\ms
{\bf Proof.} Using (2.6) and the symmetries (2.4), (2.7), we obtain
$$\om^{kl}R_{ijkl}=\om^{kl}(-R_{iklj}-R_{iljk})=\om^{lk}R_{kilj}+\om^{kl}R_{likj
}=
2K_{ij}\;.\quad \ecarre$$

{\bf Remark 3.6.} Since the left hand side of (3.5) is symmetric in $i,j$
due to
the symmetry (2.7), the formula (3.5) provides another proof of the
symmetry of the
Ricci tensor.

\ms
Yet another definition of the Ricci tensor can be extracted from the following

\ms
{\bf Proposition 3.7.} {\it Denote}
$$L=\om^{ij}\nabla_i\nabla_j\;,\leqno (3.6)$$
{\it which is considered as an operator on vector fields. Then $L$ in fact
has order $0$ and}
$$(LX)^i=L^i_{\ j}X^j,\quad L^i_{\ j}=\om^{ik}K_{kj}\;, \leqno (3.7)$$
{\it so}
$$K_{ij}=\om_{ik}L^k_{\ j}\;.\leqno (3.8)$$
{\it or on other words}
$$K(X,Y)=\om(X,LY),\quad X,Y\in T_xM\;.\leqno(3.8')$$
\ms
{\bf Proof.} We have
$$L=\om^{ij}\nabla_i\nabla_j={1\over
2}\om^{ij}(\nabla_i\nabla_j-\nabla_j\nabla_i)=
{1\over 2}\om^{ij}R(\pa_i,\pa_j)\;,$$
therefore
$$(LX)^i={1\over 2}\om^{kl}R^i_{\ pkl}X^p={1\over
2}\om^{iq}\om^{kl}R_{qpkl}X^p=
\om^{iq}K_{qp}X^p\;,$$
where we used  Proposition 3.5. Now (3.7), (3.8) and (3.8') immediately
follow. $\ecarre$
\ms
The following definition is a possibility to
introduce an analogue of the Einstein equations on a Fedosov manifold.
\bigskip
{\bf Definition 3.8.} {\it Einstein equation} on a Fedosov manifold is the
system
$$K_{ij}=0\;.\leqno(3.9)$$

\medskip\noindent
This is a system of nonlinear equations on the data of the Fedosov
manifold: its
symplectic form and connection. There is an important difference with the
symmetric
case: the system  (3.9) is not determined.
Indeed, using Darboux coordinates we can assume that $\om$ has a canonical
form and we only seek the Christoffel symbols $\Ga_{ijk}$ which should be
completely symmetric
(see Sect.1).  It is easy to see
that   the functional dimension (the number of free functional parameters)
for the space of such connections is $(n+2)(n+1)n/6$ (here $n=\dim_{\R}
M$), whereas
for the space of symmetric 2-tensors it is $(n+1)n/2$. The equality
of these functional dimensions holds if and only if $n=1$, which can not be
a dimension
of a Fedosov manifold.

\ms
On the other hand the equation (3.9) is non-trivial. Indeed, using (2.5) we
obtain
in Darboux coordinates
$$R_{likj}=\pa_k\Ga_{lji}-\pa_j\Ga_{lik}+
\om^{mp}\Ga_{pji}\Ga_{lkm}-\om^{mp}\Ga_{pik}\Ga_{ljm}\;,\leqno (3.10)$$
therefore
$$K_{ij}=\om^{kl}R_{likj}=\om^{kl}\pa_k\Ga_{lij}-\om^{kl}\om^{mp}\Ga_{pik}\Ga_{l
jm}\;,
\leqno (3.11)$$
because after multiplication by $\om^{kl}$ and summation over $k,l$ the
second and he third
terms in (3.10) vanish due to the symmetry of $\Ga_{ijk}$ in $i,j,k$.

\bigskip
{\bf Definition 3.9.} {\it Symplectic Einstein manifold} is a Fedosov manifold
such that its Ricci tensor vanishes, i.e. the Einstein equation (3.9)
is satisfied.
\ms
In terminology of I.Vaisman [Va1] such manifolds are called Ricci flat.

\ms
Let us discuss a possibility to define sectional curvature.
Consider a 2-dimensional subspace $\Pi$ in a tangent space $T_x M$ to a
Fedosov manifold $M$ at a point $x\in M$. For any $X,Y\in \Pi$ consider
a quadratic form on $\Pi$ given by
$$E_{X,Y}(Z)=R(Z,Z,X,Y)\,.\leqno(3.12)$$
Assuming that $\Pi$ is symplectically non-degenerated (non-isotropic) we can
choose $X,Y$ so that
$\omega(X,Y)=1$ and $E_{X,Y}$ has one of the following forms:
$$E_{X,Y}(Z)=r(Z_1^2+Z_2^2)\,,\quad r\in\R\setminus\{0\}\;;\leqno(3.13)$$
$$E_{X,Y}(Z)=r(Z_1^2-Z_2^2)\,,\quad r>0\;;\leqno(3.14)$$
$$E_{X,Y}(Z)=Z_1^2\,,\leqno(3.15)$$
$$E_{X,Y}(Z)=0\,,\leqno(3.16)$$
where $Z_1,Z_2$ are coordinates of $Z$ in the basis $X,Y$. Here $r=r(\Pi)$ does
not depend on the choice of the basis $X,Y$ in $\Pi$. Note that
$\omega(X,Y)=1$ implies that
the determinant of the matrix of
the form $E_{X,Y}$ in the basis $X,Y$ does not depend on $X,Y$ and equals
$r^2$ and $-r^2$ in the cases (3.13) and (3.14) respectively, and $0$ in
the cases
(3.15) and (3.16). Define also $r(\Pi)=0$ in the cases (3.15) and (3.16).

\ms
The {\it canonical form} of $E_{X,Y}$ is defined as the choice of one of
the  possibilities
(3.13)--(3.16) together with an extra parameter $r$ in cases (3.13) and (3.14)
(real non-zero or positive
respectively). We will refer to the case (3.13)  as {\it elliptic},
(3.14) as {\it hyperbolic}, (3.15) as {\it degenerate} and (3.16) as {\it
flat}.

\medskip
{\bf Definition 3.10.} {\it Sectional curvature} of a Fedosov manifold $M$ is a
function
$$\Pi\mapsto \{\hbox{canonical form of\ } E_{X,Y}\}$$
In particular we have  a continuous function $\Pi\mapsto r(\Pi)$ which is
defined
on the set of all symplectically non-degenerated
(non-isotropic)  2-dimensional subspaces $\Pi$ in the tangent bundle
$TM$.
This function may become singular
when we approach the subset of all  isotropic planes in the Grassmanian
bundle of all 2-planes
in $TM$.

\ms
{\bf Example 3.11.} Consider the 2-sphere $S^2_R$ of the radius $R$ in
$\R^3$  with the structure
of a Fedosov manifold such that
the connection is induced by the parallel transport in $\R^3$
(or, equivalently, is the Levi-Civita connection of
the Riemannian metric induced by the canonical metric in $\R^3$) and the
symplectic
form is given by the area (induced by the same Riemannian metric). Then a
straightforward
calculation shows that all tangent planes are elliptic with $r(\Pi)$
 identically equal to $1/R^2$.

\ms
{\bf Example 3.12.} Consider the standard Lobachevsky (or hyperbolic) plane
$\H^2$, e.g.
the half-plane model $\{(x,y)\in\R^2|\;y>0\}$ with the Riemannian metric
$$ds^2={dx^2+dy^2\over y^2}\;.$$
As in the previous example the corresponding area-form and Levi-Civita
connection
provide the structure of  Fedosov manifold on $\H^2$.
A straightforward
calculation shows that all tangent planes are again elliptic with $r(\Pi)=-1$
 identically.
\ms
{\bf Remark 3.13.} Assume that $M$ is a K\"ahler manifold. Then $M$ has a
canonical connection
and a canonical curvature map (2.1) (which depends on the connection only),
but it has
two curvature tensors with all the suffixes down: the Riemannian curvature
which we will denote
${\tilde R}_{ijkl}=g(\pa_i,R(\pa_k,\pa_l)\pa_j)$
and the curvature of
$M$ as a Fedosov manifold (we will denote it $R_{ijkl}$ as above). There is
an obvious
relation between them:
$${\tilde R}_{ijkl}=g_{ip}R^p_{\ jkl}=g_{ip}\om^{pm} R_{mjkl}\;. \leqno
(3.17)$$
We can write this relation in a different form using the complex structure
$J$ from (0.1)
(see [Va1]):
$${\tilde R}(X,Y,Z,W)=R(JX,Y,Z,W)\;. \leqno (3.17')$$
Indeed,
$${\tilde R}(X,Y,Z,W)=g(X,R(Z,W)Y)=\om(JX,R(Z,W)Y)=R(JX,Y,Z,W)\;.$$

\bs
\centerline{\bf 4. Normal coordinates and extensions.}

\ms
Let us consider a  manifold $M$
with a given non-degenerate 2-form $\om$.
Let $\Ga$   be a symmetric connection on $M$. (At the moment we do not
assume that
$\Ga$ preserves $\om$.) Given a point $p\in M$ we have
the exponential map $\exp_p:U\to M$ defined by $\Ga$. Here $U$ is a small
neighborhood of $0$ in $T_pM$, and $\exp_p(v)=x(1)$ where $v\in U$ and
$x(t)$ is a
geodesic defined in local coordinates $(x^1,\dots,x^n)$ near $p$ by
$${d^2x^i\over dt^2}+\Ga^i_{jk}{dx^j\over dt}{dx^k\over dt}=0,\quad x(0)=p,
\ x'(0)=v,
\leqno (4.1)$$
so $t$ is a canonical parameter along the geodesic. Denote by $(y^1,\dots,
y^n)$
the local coordinates on $M$ near $p$ which are induced by  linear coordinates
on $T_pM$ through the map $\exp_p$, so that
$$y^j(p)=0,\quad \hbox{and} \quad {\pa y^i\over \pa
x_j}(p)=\de^i_j\;.\leqno (4.2)$$
In this case $(y^1,\dots,y^n)$ are called {\it the normal coordinates
associated with
the local coordinates}
$(x^1,\dots,x^n)$ at $p$.

\ms
In other words, if $\xi=y^1{\pa\over\pa x_1}+\dots+y^n{\pa\over\pa x_n}\in
T_pM$
then $(y^1,\dots, y^n)$ are the normal coordinates of $\exp_p\xi$.

\ms
Normal coordinates are actually defined by the connection (they do not
depend on
the form $\om$). We will list several properties of the normal coordinates
(see e.g. [E], [T], [V]).

\ms
Let $(y^1,\dots, y^n)$ be local coordinates such that a given point $p\in M$
corresponds to $y^1=\dots=y^n=0$. Then they are normal coordinates if and
only if
$$\Ga^i_{jk}(y)y^jy^k\equiv 0\;,\leqno (4.3)$$
where $\Ga^i_{jk}(y)$ are the Christoffel symbols which are calculated at a
point
with the local coordinates $y=(y^1,\dots, y^n)$ in the coordinates $y^j$.

\ms
Using $\om$ we can pull the superscript $i$ down and rewrite (4.3)
equivalently in the form
$$\Ga_{ijk}(y)y^jy^k\equiv 0\;.\leqno (4.4)$$

\ms
Taking $y^j=t\xi^j$ in (4.4), dividing by $t^2$ and taking the limit as
$t\to 0$, we obtain
due to the symmetry of $\Ga_{ijk}$ in $j,k$:
$$\Ga_{ijk}(0)=0\;.\leqno (4.5)$$
\ms
Let us write the Taylor expansion of $\Ga_{ijk}$ at $y=0$ in the form:
$$\Ga_{ijk}(y)=\sum_{r=1}^\infty{1\over r!}
A_{ijk\al_1\dots \al_r}y^{\al_1}\dots y^{\al_r}\;,\leqno (4.6)$$
where
$$A_{ijk\al_1\dots\al_r}=A_{ijk\al_1\dots\al_r}(p)=
{\pa^r\Ga_{ijk}\over \pa y^{\al_1}\dots\pa y^{\al_r}}\Big|_{y=0}\;.\leqno
(4.7)$$
Here $p\in M$ is the ``origin" of the normal coordinates.

\ms

\ms
Note that different normal coordinates with the same origin differ by a linear
transformation.  It follows that the correspondence
$p\mapsto A_{ijk\al_1\dots\al_r}(p)$ defines a tensor on $M$ (the
components should be
related to our start-up coordinate system $(x^1,\dots,x^n)$).

\ms
{\bf Definition 4.1.}
The tensor $A_{ijk\al_1\dots\al_r}$ is called an  {\it affine normal
tensor} or an
{\it affine extension} of $\Ga_{ijk}$ of order $r=1,2,\dots$. We also take
by definition $A_{ijk}=0$ which is natural due to (4.5).

\bs
{\bf Proposition 4.2.} {\it The affine normal tensors have the following
properties:}

(i) $A_{ijk\al_1\dots\al_r}$ {\it is symmetric in $j,k$ and also
in $\al_1,\dots,\al_r$ (i.e. with respect to the transposition of $j$ and $k$
and also with respect to all permutations of $\al_1,\dots,\al_r$;}

(ii) {\it for any $r=1,2\dots$}
$$S_{j,k\hookrightarrow
j,k,\al_1,\dots,\al_r}(A_{ijk\al_1\dots\al_r})=0\;,\leqno (4.8)$$
{\it where $S_{j,k\hookrightarrow j,k,\al_1,\dots,\al_r}$ stands for the
sum of all
the $(r+2)(r+1)/2$ terms obtainable from the one written within parentheses
by replacing the pair $j,k$ by any  pair from the set
$\{j,k,\al_1,\dots,\al_r\}$.}

(iii) {\it Assume that a non-degenerate 2-form $\om$ on $M$ is fixed.
 If the  tensors

\ni
$\{A_{ijk\al_1\dots\al_r}|\;r=1,2,\dots\}$
are given at one point $p$,
satisfy the conditions in} (i), (ii) {\it and besides the series (4.6) is
convergent
in a neighborhood of $0$, then the sums in (4.6) constitute the set of the
Christoffel
symbols of a symmetric connection in normal coordinates. }

\ms
The proof easily follows from (4.4) and it does not differ from the proof
of the version of this statement for $\Ga^i_{jk}$
which is given e.g. in [E],[T],[V]. It can be also reduced to the
corresponding statement
for $\Ga^i_{jk}$ by pulling down the superscript $i$.

\ms
Examples of the identity (4.8) for $r=1$ and $r=2$ are as follows:
$$A_{ijkl}+A_{ijlk}+A_{iklj}=0\;;\leqno (4.9)$$
$$A_{ijklm}+A_{ijlkm}+A_{ijmkl}+A_{ikljm}+A_{ikmjl}+A_{ilmjk}=0\;.
\leqno (4.10)$$

\ms
Let us consider an arbitrary tensor $T=(T^{i_1\dots i_k}_{j_1\dots j_l})$
on $M$.

\ms
{\bf Definition 4.3.} The (affine) {\it extension} of $T$ of order $r$ is a
tensor on $M$,
such that its components at the point $p\in M$ in the local coordinates
$(x^1,\dots,x^n)$
are given by the formulas
$$T^{i_1\dots i_k}_{j_1\dots j_l,\al_1\dots\al_r}=
{\pa^rT^{i_1\dots i_k}_{j_1\dots j_l}
\over \pa y^{\al_1}\dots\pa y^{\al_r}}\Big|_{y=0}\;,\leqno (4.11)$$
where $(y^1,\dots,y^n)$ are normal coordinates associated with
$(x^1,\dots,x^n)$ at $p$.

\ms
In particular we can form the  extensions of $\om_{ij}$ which will be tensors
$\om_{ij,\al_1\dots\al_r}$.

\ms
Note that the first extension of any tensor coincides with its covariant
derivative
because $\Ga^i_{jk}(0)=0$ in normal coordinates.

\ms
Let us consider the extensions of $\om$ and $\Ga$ together.

\ms
{\bf Lemma 4.4.} (i) {\it If the connection $\Ga$ preserves $\om$ then}
$$\om_{ij,k\al_1\dots\al_r}=A_{ikj\al_1\dots\al_r}-A_{jki\al_1\dots\al_r}
\leqno (4.12)$$
{\it for all $r=0,1,\dots$ and all $i,j,k,\al_1,\dots,\al_r$. }

(ii) {\it Vice versa if (4.12) holds  on $M$ for $r=0$  then $\Ga$
preserves $\om$.}

(iii) {\it If $M,\om,\Ga$ are real analytic and the relations (4.12) hold
at a point $p\in M$
then $\Ga$ preserves $\om$ in the connected component of $p$.}

\ms
{\bf Proof.} (i) Let us recall (see Sect.1) that $\Ga$ preserves $\om$ if
and only if
$$\pa_k\om_{ij}=\Ga_{ikj}-\Ga_{jki}\;.\leqno (4.13)$$

Let $(y^1,\dots,y^n)$ be normal coordinates at $p\in M$. Let us write the
Taylor
expansion of $\om_{ij}$ in these coordinates
$$\om_{ij}(y)=\sum_{r=0}^\infty {1\over
r!}\om_{ij,\al_1\dots\al_r}(0)y^{\al_1}\dots y^{\al_r}\;.
\leqno (4.14)$$
It follows that $\pa_k\om_{ij}={\pa\over\pa y^k}\om_{ij}$ has the Taylor
expansion
$$\pa_k\om_{ij}(y)=\sum_{r=1}^\infty {1\over r!}\om_{ij,\al_1\dots\al_r}(0)
\sum_{s=1}^r y^{\al_1}\dots \de_k^{\al_s}\dots y^{\al_r}$$
$$=
\sum_{r=1}^\infty {1\over
r!}\sum_{s=1}^r\om_{ij,k\al_1\dots\hat\al_s\dots\al_r}(0)
y^{\al_1}\dots \hat y^{\al_s}\dots y^{\al_r}\;.$$
(Here a ``hat" \^\ over a factor means that this factor should be omitted.)
Clearly all the terms in the last sum are equal, therefore we obtain
$$\pa_k\om_{ij}(y)=\sum_{r=1}^\infty {1\over
(r-1)!}\om_{ij,k\al_1\dots\al_{r-1}}(0)
y^{\al_1}\dots \dots y^{\al_{r-1}}$$
or
$$\pa_k\om_{ij}(y)=\sum_{r=0}^\infty {1\over r!}\om_{ij,k\al_1\dots\al_r}(0)
y^{\al_1}\dots \dots y^{\al_r}\;.\leqno (4.15)$$
Now comparing the Taylor decompositions in both sides of (4.13) we arrive
to (4.12).

\ms
(ii) Assuming that (4.12) holds for $r=0$, we see that the first extension
of $\om_{ij}$
vanishes on $M$, and the statement (ii) follows because this extension
coincides
with the covariant derivative.

\ms
(iii) If $M,\om,\Ga$ are real analytic then the series (4.14) converges and
can be
termwise differentiated which leads to (4.15) (where the series is also
convergent).
Comparing the Taylor series of the both sides of (4.13) we see that the
relations (4.12)
at $p$ imply that (4.13) holds in a neighborhood of $p$ which proves (iii).
$\ecarre$

\ms
Now we can formulate a point characterization of  possible extensions of
$\om_{ij}$
and $\Ga_{ijk}$ at a point of a Fedosov manifold.

\ms
{\bf Theorem 4.5.} 1) {\it Assume that we are given a set of tensors }
$$\om_{ij,\al_1\dots\al_r}=\om_{ij,\al_1\dots\al_r}(0),\
A_{ijk\al_1\dots\al_r}=A_{ijk\al_1\dots\al_r}(0);\ r=0,1,\dots \leqno (4.16)$$
{\it at a point $0\in\R^n$. Then a structure of a real-analytic Fedosov
manifold
in a neighborhood of $0$ in $\R^n$ with the tensors (4.16) serving as
extensions of
$\om_{ij}$ and $\Ga_{ijk}$ in normal coordinates at $0$, exists if and only
if the following
conditions are satisfied:}
\ms
(a) {\it $A_{ijk}(0)=0$, and $A_{ijk\al_1\dots\al_r}(0)$ is symmetric
in $j,k$ as well as in $\al_1,\dots,\al_r$ for any $r=1,2\dots$;}

(b) $A_{ijk\al_1\dots\al_r}(0)$ {\it satisfy the identity (4.8);}

(c) $\om_{ij}(0)$ {\it is skew-symmetric in $i,j$ and non-degenerate;}

(d) $\om_{ij,\al_1\dots\al_r}(0)$ {\it is skew-symmetric in $i,j$ and
symmetric in $\al_1,\dots,\al_r$;}

(e) {\it for all $r=0,1,\dots$ the identity (4.12) is satisfied at $0$;}

(f) {\it the series (4.6) and (4.14) are convergent for small $y$.}

\ssni
{\it The  structure of a real-analytic Fedosov manifold near $0$ is
uniquely defined
by the tensors (4.16).}
\ms
2) {\it Assume that we are only given tensors $A_{ijk\al_1\dots\al_r}(0)$
such that the conditions (a), (b) are satisfied and the series (4.6) is
convergent for
small $y$. Then a real-analytic symplectic form near $0$ complementing the
connection
(defined by $\Ga_{ijk}$ given by the series (4.6)) to a local real-analytic
Fedosov manifold
structure near $0$ exists if and only if for all $r=1,2,\dots$
$$A_{ikj\al_1\dots\al_r}(0)-A_{jki\al_1\dots\al_r}(0)\ \hbox{is symmetric in}
\ k,\al_1\;,\leqno (4.18)$$
i.e. with respect to the transposition of $k$ and $\al_1$.
The form $\om$ above is uniquely defined by $\om_{ij}(0)$ and it exists for
arbitrary
skew-symmetric non-degenerate $\om_{ij}(0)$.
}

\ms
{\bf Proof.} 1) It is obvious from the considerations above that the
conditions (a)-(f)
are necessary for the existence of a local real-analytic Fedosov manifold
near $0$
with the given extensions (4.16). They are also sufficient because we can
define
$\om_{ij}(y)$ and $\Ga_{ijk}(y)$ near $0$ by the convergent series (4.6)
and (4.14),
and it is easy to see that all the conditions will be  satisfied.
The uniqueness of $\om,\Ga$ near
$0$ is obvious because the tensors (4.16) define all the Taylor
coefficients of $\om,\Ga$
at $0$.

\ms
2) If the tensors $A_{ijk\al_1\dots\al_r}$ satisfy (a),(b) and the series
(4.6) are
convergent for small $y$, then the sums of these series provide $\Ga_{ijk}(y)$
which determine a real-analytic symmetric connection in normal coordinates
for any given real-analytic non-degenerate skew-symmetric form
$\om_{ij}(y)$ defined near $0$. To make this 2-form symplectic and
preserved by $\Ga$
we can try to find it from the conditions (4.13) which determine all first
derivatives $\pa_k\om_{ij}$ near $0$, so we only need to
know $\om_{ij}(0)$.

Let us choose an arbitrary skew-symmetric matrix $\om_{ij}(0)$.
Then we can find $\om_{ij}(y)$ from (4.13) if and only if the following
compatibility
conditions are satisfied near $0$:
$$\pa_l(\Ga_{ikj}-\Ga_{jki})=\pa_k(\Ga_{ilj}-\Ga_{jli})\;,\leqno (4.19)$$
or, equivalently,
$$\pa_l(\Ga_{ikj}-\Ga_{jki})\ \hbox{is symmetric in}\ k,l\;.\leqno (4.20)$$
But the same  argument as the one leading to (4.15), shows that the Taylor
expansion
of $\pa_l\Ga_{ikj}(y)$ at $0$ has the form
$$\pa_l\Ga_{ikj}(y)=
\sum_{r=1}^\infty {1\over r!}A_{ikjl\al_1\dots\al_r}y^{\al_1}\dots y^{\al_r}\;.
\leqno (4.21)$$
It follows that (4.20) is equivalent to (4.18), so the statement 2)
follows. $\ecarre$

\ms
{\bf Corollary 4.6.} {\it Let $A_{ijk\al_1\dots\al_r}$ be tensors at a
point $0\in \R^n$,
$r=0,1,\dots$, $A_{ijk}=0$ and the series (4.6) are convergent near $0$.
Then the sum of these series define symplectic connection in normal
coordinates,
which are at the same time Darboux coordinates
(with arbitrary constant skew-symmetric non-degenerate $\om_{ij}$), if and
only if
$A_{ijk\al_1\dots\al_r}$ are symmetric in $i,j,k$ and in $\al_1,\dots,\al_r$,
and besides the conditions (4.8) are satisfied.}

\ms
{\bf Proof.} Let us fix an arbitrary constant skew-symmetric non-degenerate
$\om_{ij}$.
The condition (4.8) guarantees that the sums $\Ga_{ijk}(y)$ of the series
(4.6) give us
components of a symmetric connection in normal coordinates due to
Proposition 4.2.
Now (1.2) implies that this connection preserves the constant form $\om$ if
and only if
$\Ga_{ijk}$ is symmetric in $ijk$ which is equivalent to saying that
$A_{ijk\al_1\dots\al_r}$ are symmetric in $i,j,k$ for any $r=1,2,\dots\;.$
$\ecarre$

\ms
Now we can formulate the simplest result about local invariants of Fedosov
manifolds.
For simplicity we will only consider invariants which are tensors though
similar
results are true e.g. for invariants which are tensors with values in
density bundles.

\ms
{\bf Definition 4.7.} A tensor $T=(T^{i_1\dots i_k}_{j_1\dots j_l}(x))$ on
a Fedosov manifold
$M$ is called a {\it local invariant} of $M$ if its components at each
point $p\in M$
are functions of a finite
number of derivatives (of arbitrary order, including $0$) of $\Ga_{ijk}$
and $\om_{ij}$
taken at the same point $p$ in local coordinates, i.e.
$$T^{i_1\dots i_k}_{j_1\dots j_l}(p)=F^{i_1\dots i_k}_{j_1\dots j_l}
(\Ga_{ijk}(p),\pa_\al\Ga_{ijk}(p),
\dots,\om_{ij}(p),\pa_\al\om_{ij}(p),\pa_{\al_2}\pa_{\al_1}\om_{ij}(p),\dots)\;,
\leqno (4.22)$$
where the functions $F^{i_1\dots i_k}_{j_1\dots j_l}$ do not depend on the
choice of the local
coordinates. A local invariant is called {\it rational} or {\it polynomial}
if the functions
$F^{i_1\dots i_k}_{j_1\dots j_l}$ are all rational or polynomial  respectively.

\ms
In particular, this definition applies to local invariants which are
differential forms.
\ms
{\bf Theorem 4.8.} {\it Any local invariant of a Fedosov manifold is a function
 of $\om_{ij}$ and  a finite number of its affine normal tensors, i.e.
it can be presented in the form}
$$T^{i_1\dots i_k}_{j_1\dots j_l}(p)=G^{i_1\dots i_k}_{j_1\dots j_l}
(\om_{ij}(p), A_{ijk\al_1}(p),\dots,A_{ijk\al_1\dots\al_r}(p))\;,\leqno
(4.23)$$
{\it where the functions $G^{i_1\dots i_k}_{j_1\dots j_l}$ do not depend on
the choice
of local coordinates. If the invariant $T^{i_1\dots i_k}_{j_1\dots j_l}$ is
rational
or polynomial, then the functions $G^{i_1\dots i_k}_{j_1\dots j_l}$ are
also rational
or polynomial respectively.}

\ms
{\bf Proof.} Since $T$ is a tensor, its components $T^{i_1\dots
i_k}_{j_1\dots j_l}(p)$
will not change if we replace local coordinates $(x^1,\dots, x^n)$ by
normal coordinates
$(y^1,\dots, y^n)$ which are associated with $(x^1,\dots, x^n)$ at $p$. But
then according
to  Definition 4.7 we obtain, using formulas (4.7) and (4.12):
$$T^{i_1\dots i_k}_{j_1\dots j_l}(p)=
F^{i_1\dots i_k}_{j_1\dots j_l}(0,A_{ijk\al}(p),\dots,
\om_{ij}(p),0,A_{i\al_1j\al_2}(p)-A_{j\al_1i\al_2}(p),\dots)\;.
\leqno (4.24)$$
All statements of the Theorem immediately follow. $\ecarre$

\bs
\centerline {\bf 5. Normal tensors and  curvature.}
\bs
{\bf A.} Let us try to connect the affine  normal tensors and extensions of
$\om_{ij}$
on a Fedosov
manifold $M$ with the curvature tensor. First let us rewrite (2.5) in terms of
$\Ga_{ijk}$ and $R_{ijkl}$ instead of
$\Ga^i_{jk}$ and $R^i_{\ jkl}$:
$$
R_{ijkl}=\om_{is}\pa_k(\om^{sp}\Ga_{pjl})-\om_{is}\pa_l(\om^{sp}\Ga_{pjk})-
\om^{mp}\Ga_{pjl}\Ga_{ikm}+\om^{mp}\Ga_{pjk}\Ga_{ilm}\;. \leqno (5.1)
$$
The standard formula for differentiation of the inverse matrix gives
$$\pa_k\om^{ij}=-\om^{ir}\om^{sj}\pa_k\om_{rs}=\om^{ir}\om^{js}\pa_k\om_{rs}\;.
\leqno (5.2)$$

Let us introduce normal coordinates with the center $p$. Observe that
in such coordinates $\Ga_{ijk}(p)=(\pa_k\om_{ij})(p)=0$. It  follows from
(5.2)
that $(\pa_k\om^{ij})(p)=0$. Now performing differentiation in the right
hand side
of (5.1) and taking the values at $p$ we obtain
$$R_{ijkl}=A_{ijlk}-A_{ijkl}\;,\leqno (5.3)$$
which is an equality of tensors, and it is true everywhere on $M$ because
$p\in M$ is
an arbitrary point. This gives an expression of $R_{ijkl}$ as a
linear combination of affine normal tensors.

\ms
We can do another step in this direction by applying $\pa_m$ to both sides
of (5.1)
and taking values at $p$ which leads to
$$R_{ijkl,m}=A_{ijlkm}-A_{ijklm}\;,\leqno (5.4)$$
where $R_{ijkl,m}$ means the first extension of the curvature tensor
$R_{ijkl}$ or
its first covariant derivative (which coincide).

\ms
Further differentiations lead to expressions of the extensions of the
curvature tensor
which are already non-linear in normal tensors and also explicitly contain
$\om^{ij}$,
though they are linear in the highest order normal tensors:

\ms
{\bf Lemma 5.1.} {\it For any $r\ge 0$}
$$
R_{ijkl,\al_1\dots \al_r}=
A_{ijlk\al_1\dots \al_r}-A_{ijkl\al_1\dots \al_r}+P_r \leqno (5.5)
$$
{\it where $P_r$ is a polynomial of the components $\om^{ij}$ and of the
normal tensors of order $\le r$, $P_0=P_1=0$.}

\ms
{\bf B.} Our goal now is to express the normal tensors and the extensions
of $\om_{ij}$
in terms of the curvature tensor. Let us start by solving the equations (5.3)
with respect to $A_{ijkl}$. Making cyclic
permutations of $j,k,l$ we obtain, using  (4.9):
$$R_{iklj}=A_{ikjl}-A_{iklj}=A_{ijkl}-A_{iklj}=2A_{ijkl}+A_{ijlk}\;,\leqno
(5.3')$$
$$R_{iljk}=A_{ilkj}-A_{iljk}=A_{iklj}-A_{ijlk}=-2A_{ijlk}-A_{ijkl}\;.\leqno
(5.3'')$$
Adding doubled (5.3') with (5.3'') gives
$$A_{ijkl}={1\over 3}(2R_{iklj}+R_{iljk})\;.\leqno (5.6)$$
Also using (2.6), we obtain
$$A_{ijkl}={1\over 3}(R_{iklj}+R_{ijlk})\;.\leqno (5.7)$$

A simple relation between the second extension of $\om_{ij}$ and the
curvature tensor is
given by

\ms
{\bf Proposition 5.2.} {\it On any Fedosov manifold}
$$\om_{ij,kl}={1\over 3}R_{klij}\;.\leqno (5.8)$$

\ms
{\bf Proof.} According to Lemma 4.4 we obtain, using (5.7):
$$\om_{ij,kl}=A_{ikjl}-A_{jkil}={1\over
3}(R_{iklj}+R_{ijlk}-R_{jkli}-R_{jilk})\;.$$
The second and  fourth terms cancel due to Proposition 2.3, so using (2.6)
we obtain:
$$\om_{ij,kl}={1\over 3}(R_{iklj}-R_{jkli})={1\over 3}(R_{kilj}-R_{kjli})=
{1\over 3}(-R_{kijl}-R_{kjli})={1\over 3}R_{klij}\;.\quad \ecarre$$

Now we will try to solve the equations (5.5) with respect to the higher order
normal tensors.
 Relations of this type between
the curvature tensor $R^i_{\ jkl}$ and normal tensors $A^i_{jk\al_1\dots\al_r}$
defined similarly to (4.7) but with $\Ga_{ijk}$ replaced by $\Ga^i_{jk}$,
can be found in
[Th], Sect.49, and [V], Ch. VI, Sect.9-12. Note that generally
$A^i_{jk\al_1\dots\al_r}$
does not coincide with
$\om^{ip}A_{pjk\al_1\dots\al_r}$  because $\om^{ip}$ is not constant.
We will give more details than in the exposition in [Th] and [V] which we
found
insufficient.

\ms
We will need a combinatorial preparation first.
Let $J=\{j_1j_2\dots ,j_m\}$ be an ordered set of positive
integers. We associate to $J$ an ordered set
$T(J)$ of ordered triads of elements of $J$. The set $T(J)$
will be defined by an inductive procedure:

$$i)\ T(J)=\emptyset,\ m=1,2; $$

$$ii)\ T(J)=\{(j_1j_2j_3)<(j_3j_1j_2)\},\ m=3; $$

iii) For $m>3$ the elements $T(J)$ in the increasing order are:
$$
(j_1j_2j_3), (j_1j_3j_4),\dots , (j_1j_{m-1}j_m),
$$
$$(j_mj_1j_2),$$
$$
(j_2j_mj_{m-1}), (j_2j_{m-1}j_{m-2}),\dots , (j_2j_4j_3),
$$
$$(j_3j_2j_4), $$
$$T(J\setminus \{j_1j_2\}).$$

\ms
{\bf Remark 5.3.} The most important property of this ordering is the
following:
the first and third elements of any triad are the same pair (with possibly
reversed order)
as the first two elements of the next triad.

\ms
{\bf Lemma 5.4.} {\it For $m\ge 2$ we have}
$$\hbox{Card}(T(J))=m(m-1)/2 - 1\;.\leqno (5.9)$$

\ms
{\bf Proof}. The formula holds for $m=2$ and $m=3$. Now let us use induction
with respect to $m$. Denote
$N_m=\hbox{Card}(T(M))$. By the construction

$$
N_m=(m-2)+1+(m-3)+1+N_{m-2}=N_{m-2}+(2m-3).
$$
On the other hand it is easy to check that

$$m(m-1)/2 = (m-2)(m-3)/2 +(2m-3).\qquad \ecarre$$

\ms
Let  $R_{ijkl,\al_1\dots \al_r}$ be the $r$-th
extension of the curvature tensor $R_{ijkl}$ on a Fedosov manifold.
Consider the ordered set $J=\{jkl\al_1\dots \al_r\}$.
Note that we ignore the first subscript $i$ which plays a special role.
For any triad $u=(u_1u_2u_3)\in T(J)$ consider the
$r$-th extension $\tilde R_{iu}$ of the curvature tensor
$R_{iu}$:

$$\tilde R_{iu}=R_{iu, J\setminus u}.$$

Denote $N=\hbox{Card}(T(J))$ and let $u^{(1)}<u^{(2)}<\dots <u^{(N)}$
be the list of all elements of $T(J)$ (in increasing order).

\ms
The following theorem is an analogue of a Veblen theorem from [V].

\ms
{\bf Theorem 5.5.} {\it For any $r\ge 0$ the affine normal tensor of order
$r+1$
can be expressed through the extensions of the curvature tensor as follows:}

$$
A_{ijkl\al_1\dots \al_r}=
-{1\over N+1}\sum_{s=1}^N (N-s+1)\tilde R_{iu^{(s)}}+P_r\;,
\leqno (5.10)
$$
{\it where $P_r$ is a polynomial of $\om^{ij}$ and affine normal tensors of
order
$p\le r$. Besides

\ni
$P_0=P_1=0$, so for $r=0,1$ we have precise formulas}
$$A_{ijkl}=-{1\over 3}(2R_{ijkl}+R_{iljk})\;; \leqno (5.10)_0$$
$$A_{ijklm}=-{1\over
6}(5R_{ijkl,m}+4R_{ijlm,k}+3R_{imjk,l}+2R_{ikml,j}+R_{ilkm,j})\;.
\leqno (5.10)_1$$

\ms
To prove the theorem we will need  several simple facts.
We recall, first of all, that the normal tensors
$A_{ijk\al_1,\dots ,\al_r}$ are symmetric in
$j,k$ and also in  $\al_1,\dots ,\al_n$.

Let us fix the ordered set of indices $J=\{j,k,\al_1\dots ,\al_r\}$.
Denote by $P(J)$ the set of all ordered pairs $p=(p_1,p_2)$
of elements from $J$ (here ``ordered" means that $p_1$ stands somewhere
before $p_2$
in the list of the elements of $J$ above). Denote by $\tilde A_{ip}$ the
normal tensor $A_{ip_1p_2,J\setminus p}$ of degree $r$. The equality (4.8)
(in Proposition 4.2) can be rewritten as follows:

$$\sum_{p\in P(J)}\tilde A_{ip} =0.\leqno (5.11)$$

\ms
{\bf Lemma 5.6.} (i) {\it The extension of order $r\ge 1$ of the tensor
$\om_{ij}$
is a linear combination  of normal tensors of order $r-1$. }

(ii) {\it Every component of the extension of order $r\ge 1$ of the tensor
$\om^{ij}$
is a polynomial (independent on the choice of the coordinates)
 of the components $\om^{ij}$ and normal tensors of order $\le r-1$.}

\ms
{\bf Proof.} The  part (i) is obvious from Lemma 4.4 (in fact the formula
(4.12)
gives an explicit expression).  The part (ii) is obtained by multiple
differentiation
of the  formula (5.2)
in normal coordinates with  substitution of the right hand side
of (5.2) instead of the first derivatives of $\om^{ij}$ after each
differentiation,
and further use of the part (i) (or formula (4.12)). $\ecarre$

\ms
{\bf Proof of Theorem 5.5.} Let
$J=\{j,k,l,\al_1,\dots , \al_r\}$.
Applying Lemma 5.1 and using
the symmetries of normal tensors
one could write that

$$
-\sum_{s=1}^N (N-s+1)\tilde R_{iu^{(s)}}= N\tilde A_{ijk}-
$$
$$
- (N\tilde A_{ijl} - (N-1)\tilde A_{ijl}) -
$$
$$
- ((N-1)\tilde A_{ij\al_1}- (N-2)\tilde A_{ij\al_1})- \dots =
$$
$$
=N\tilde A_{ijk} -
\sum_{p\in P(J),\ p\neq (jk)}\tilde A_{ip}
$$
modulo a polynomial of $\om^{ij}$ and normal tensors of order $\le r$
(this polynomial vanishes if $r=0$ or 1).
Now the theorem follows because the right hand side equals $(N+1)\tilde
A_{ijk}$
due to (5.11).
$\ecarre $

\ms
{\bf Corollary 5.7.} {\it Any  component $A_{ijk\al_1\dots\al_r}$ of the
normal tensor
of order $r$ can be expressed
as a universal polynomial (i.e. a polynomial independent of the coordinates)
of  components $\om^{ij}$ and components of extensions of the curvature tensor
$R_{ijkl}$ of order $\le r-1$.}

\ms
{\bf Proof.} The result follows from Theorem 5.5 if we use (5.6) (or (5.7)) and
then implement induction in $r$. $\ecarre$

\ms
{\bf C.} Now we will establish relations between extensions and covariant
derivatives.
Let us consider an arbitrary tensor $T=(T^{i_1\dots i_k}_{j_1\dots j_l})$.
Let us recall
that its (first order) covariant derivative is a tensor
$\nabla_{\pa_\al}T^{i_1\dots i_k}_{j_1\dots j_l}=
\nabla_\al T^{i_1\dots i_k}_{j_1\dots j_l}=
T^{i_1\dots i_k}_{j_1\dots j_l,\al}$
defined by the formulas
$$
T^{i_1\dots i_k}_{j_1\dots j_l,\al}=
\nabla_\al T^{i_1\dots i_k}_{j_1\dots j_l}=
\nabla_{\pa_\al}T^{i_1\dots i_k}_{j_1\dots j_l}=
 \leqno (5.12)
$$

\ni
$
\pa_\al T^{i_1\dots i_k}_{j_1\dots j_l}+
\Ga^{i_1}_{\al s}T^{si_2\dots i_k}_{j_1\dots j_l}+
\dots+\Ga^{i_k}_{\al s}T^{i_1\dots i_{k-1}s}_{j_1\dots j_l}-
\Ga^s_{\al j_1}T^{i_1\dots i_k}_{sj_2\dots j_l}-\dots -
\Ga^s_{\al j_l}T^{i_1\dots i_k}_{j_1\dots j_{l-1}s}=$

\msni
$
\pa_\al T^{i_1\dots i_k}_{j_1\dots j_l}+
\om^{i_1q}\Ga_{q\al s}T^{si_2\dots i_k}_{j_1\dots j_l}+
\dots+\om^{i_kq}\Ga_{q\al s}T^{i_1\dots i_{k-1}s}_{j_1\dots j_l}-
$
$$
\om^{sq}\Ga_{q\al j_1}T^{i_1\dots i_k}_{sj_2\dots j_l}-\dots-
\om^{sq}\Ga_{q\al j_l}T^{i_1\dots i_k}_{j_1\dots j_{l-1}s}\;.
$$
In particular, taking normal coordinates we immediately see that the covariant
derivative of $T$ coincides with the first extension of $T$.

We  define
the higher order covariant derivative of $T$ of order $r$ as the following tensor:
$$T^{i_1\dots i_k}_{j_1\dots j_l,\al_1,\al_2,\dots,\al_r}=\nabla_{\al_r
}\dots\nabla_{\al_1}T^{i_1\dots i_k}_{j_1\dots j_l}\;.\leqno (5.13)$$
It depends on the order of $\al_1,\dots,\al_r$ (unlike the extension
$T^{i_1\dots i_k}_{j_1\dots j_l,\al_1\dots\al_r}$).

\ms
{\bf Proposition 5.8.} (i) {\it Every component of the covariant derivative
of $T$
of order $r$ can be expressed through extensions of $T$ of order $\le r$:}
$$T^{i_1\dots i_k}_{j_1\dots j_l,\al_1,\al_2,\dots,\al_r}=
T^{i_1\dots i_k}_{j_1\dots j_l,\al_1\al_2\dots\al_r}+P_{r-1}\;,\leqno (5.14)$$
{\it where $P_{r-1}$ is a linear combination of components of extensions of $T$
of order $\le r-1$ with coefficients which are universal polynomials of
$\om^{ij}$ and of
the components of the normal tensors of order $\le r-1$.}

(ii) {\it Every component of the extension of $T$ of order $r$ can be
expressed through
the covariant derivatives of $T$of order $\le r$:}
$$T^{i_1\dots i_k}_{j_1\dots j_l,\al_1\al_2\dots\al_r}=
T^{i_1\dots i_k}_{j_1\dots j_l,\al_1,\al_2,\dots,\al_r}+P_{r-1}\;,\leqno
(5.15)$$
{\it where $P_{r-1}$ is a linear combination of components of covariant
derivatives of $T$
of order $\le r-1$ with coefficients which are universal polynomials of
$\om^{ij}$ and of
the components of the normal tensors of order $\le r-1$.}

\ms
{\bf Proof.} (i) is obtained by repeated application of covariant
differentiation
formula (5.12) in normal coordinates with the use of (5.2) to differentiate
$\om^{ij}$.
We should also use (1.2) to get rid of the first derivatives of $\om_{ij}$
as soon as they show up, replacing them by linear combinations of $\Ga_{ijk}$.
After all differentiations are done, we should take values of all
quantities at the
center of the normal coordinates, which leads to (5.14).

\ms
To prove (ii) we should use induction in $r$. We  already observed that
(5.15) is true
for $r=1$ (with $P_0=0$), i.e. the first extension coincides with the first
covariant derivative. Now assuming that (5.15) is true with $r\le r_0-1$,
we can take
(5.14) with $r=r_0$ and replace components of the  extensions of $T$ in
$P_{r-1}$
by the corresponding combinations of the components of the covariant
derivatives.
This leads to (5.15) with $r=r_0$. $\ecarre$

\ms
{\bf Theorem 5.9.} {\it Any component $A_{ijk\al_1\dots\al_r}$ of the
normal tensor
of order $r$ can be expressed
as a universal polynomial (i.e. a polynomial independent of the coordinates)
of the components $\om^{ij}$ and of the components of covariant derivatives
of the curvature tensor $R_{ijkl}$ of order $\le r-1$.}

\ms
{\bf Proof.} Let us act by induction in $r$. For $r=1$ the statement
follows from
(5.6) (or (5.7)). For any $r\ge 2$ let us use Theorem 5.5 to express any
component
of the normal tensor
of order $r+1$ through  extensions of the curvature tensor of order $\le r$
by the formula (5.10). In the remainder term $P_r$ of (5.10) we should
replace the normal
tensors of order $\le r$ by the covariant derivatives of the curvature tensor
of order $\le r-1$ using the induction assumption. Now express the extensions
$\tilde R_{iu^{(s)}}$  in the right hand side of (5.10) by the covariant
derivatives
of the curvature tensor using Proposition 5.8.  It remains to use again the
induction
assumption to replace
the newly appeared components of the normal tensors (from the remainders in
(5.15))
by the components of the covariant derivatives of $R_{ijkl}$. $\ecarre$

\ms
{\bf Corollary 5.10.} {\it Proposition 5.8. remains true if in its
statement we replace
the components of the normal tensors of order $\le r-1$ by the components
of the covariant
derivatives of order $\le r-2$ of the curvature tensor $R_{ijkl}$. }

\ms
Now we are ready for another presentation of possible local invariants of a
Fedosov manifold.
\ms
{\bf Theorem 5.11.} {\it Any local invariant of a Fedosov manifold is a
function
 of $\om_{ij}$ and  a finite number of covariant derivatives of its
curvature tensor
$R_{ijkl}$, i.e.
it can be presented in the form}
$$T^{i_1\dots i_k}_{j_1\dots j_l}(p)=G^{i_1\dots i_k}_{j_1\dots j_l}
(\om_{ij}(p), R_{ijkl}(p),\dots,R_{ijkl,\al_1,\dots,\al_r}(p))\;,\leqno
(5.16)$$
{\it where the functions $G^{i_1\dots i_k}_{j_1\dots j_l}$ do not depend on
the choice
of local coordinates. If the invariant $T^{i_1\dots i_k}_{j_1\dots j_l}$ is
rational
or polynomial (in $\om_{ij},\om^{ij}, \Ga_{ijk}$ and derivatives of
$\om_{ij}$, $\Ga_{ijk}$),
then the functions $G^{i_1\dots i_k}_{j_1\dots j_l}$ are also rational
or polynomial (in the same sense) respectively.}

\ms
{\bf Proof.} Due to Theorem 4.8 and Proposition 5.2 it is sufficient to
express the normal
tensors $A_{ijk\al_1\dots\al_r}$ polynomially in terms of $\om^{ij}$
and the covariant derivatives  $R_{ijkl,\al_1,\dots,\al_p}$, $p=0,\dots, r-1$.
without $\om$)
induction
This is possible due to Theorem 5.9.
$\ecarre$

\ms
{\bf Remark 5.12.} Not all functions (and even not all polynomials)
$G^{i_1\dots i_k}_{j_1\dots j_l}$ indeed give us local invariants of the
Fedosov
manifold (the corresponding expression is not necessarily a tensor yet, it
might still
depend on the choice of local coordinates). The form given in the
expression (5.16)
is necessary but not sufficient for an expression to be a local invariant.
We will describe sufficient conditions later in this paper.

\ms
{\bf D.} Theorem 4.5 establishes necessary and sufficient conditions for
the tensors

\ni
$A_{ijk\al_1\dots\al_r}$, given at a point, to be normal tensors of a
real-analytic
Fedosov manifold. The relations between normal tensors and the extensions
of the
curvature tensor or with its covariant derivatives which we proved above,
in principle
allow to write necessary and sufficient conditions for a family of tensors
to be
a family of extensions or covariant derivatives of the curvature tensor of
a real-analytic Fedosov manifold at a point.  However the algebraic
relations between
normal tensors and the curvature are not simple, so we can do it effectively
on the lowest levels only. As an example we will do it for the curvature tensor
itself (without extensions or derivatives) and also for the first covariant
derivatives
of the curvature tensor.

\ms
{\bf Theorem 5.13.} {\it Assume that we are given a tensor $R_{ijkl}$ at
the origin
$0\in\R^n$. Then  necessary and sufficient conditions for this tensor to be
the curvature
tensor of  a Fedosov manifold at a point are:}

$$ R_{ijkl}=-R_{ijlk}\;; \leqno \hbox{(a)}$$
$$ R_{ijkl}=R_{jikl}\;; \leqno \hbox{(b)}$$
$$R_{i(jkl)}:= R_{ijkl}+R_{iklj}+R_{iljk}=0\;.\leqno \hbox{(c)}$$

\ms
{\bf Proof.} The conditions (a), (b) (c) are necessary as was established
in Sect. 2.
Let us prove that they are sufficient.

Let us introduce a  hypothetical normal
tensor $A_{ijkl}$ at $0$ by the formulas (5.6) or (5.7) (it follows from
(a) and (c) that
they are equivalent). Let us check that the tensor $A_{ijkl}$ satisfies the
conditions
from the part 2) of Theorem 4.5 with with $r=1$ i.e.
$$A_{ijkl}=A_{ikjl} \;;\leqno (5.17)$$
$$A_{i(jkl)}=0\;;\leqno (5.18)$$
$$A_{ijkl}-A_{kjil}-A_{ilkj}+A_{klij}=0\;.\leqno (5.19)$$
Clearly (5.17) follows from (5.7), and (5.18) follows from the condition
(c). Now
substituting the expressions for the terms in the left hand side of (5.19)
from (5.7)
we get
$${1\over 3}(R_{iklj}+R_{ijlk}-R_{kilj}-R_{kjli}-R_{ikjl}-R_{iljk}+
R_{kijl}+R_{klji})=-{1\over 3}(R_{ijkl}+R_{jkli}+R_{klij}+R_{lijk}),$$
where we used the properties (a) and (b). The last expression vanishes due to
Proposition 2.4 (it is established in its proof that (2.9) follows from
the properties (a), (b), (c)).

Now defining $A_{ijk\al_1\al_2\dots\al_r}=0$ for all $r\ge 2$ we obtain a
system of
tensors which can serve as a system of normal tensors of a Fedosov manifold
at a point
due to Theorem 4.5. It remains to prove that its curvature tensor $\tilde
R_{ijkl}$
at the given point will
coincide with the given tensor $R_{ijkl}$. By (5.2) we have
$$\tilde R_{ijkl}=A_{ijlk}-A_{ijkl}={1\over
3}(2R_{ilkj}+R_{ikjl}-2R_{iklj}-R_{iljk})=
{1\over 3}(3R_{ilkj}+3R_{ikjl})=R_{ijkl}\;,$$
where we used  the properties (a) and (c). $\ecarre$

\ms
To work with the first derivatives we will need another expression of the
normal tensors
$A_{ijklm}$ through the first covariant derivatives (or first extensions)
of $R_{ijkl}$.
It is similar to $(5.10)_1$ but has smaller coefficients.

\ms
{\bf Lemma 5.14.} {\it On any Fedosov manifold}
$$A_{ijklm}=-{1\over
6}(2R_{ijkl,m}+R_{ijkm,l}+R_{ikjl,m}+R_{ikjm,l}+R_{iljm,k}).
\leqno (5.20)$$

\ms
{\bf Proof.} The formula (5.20) can be deduced from a similar formula
(49.7) in [Th], p.130,
by pulling down the superscript. However for the convenience of the reader
we will
supply the proof which is only sketched in [Th].

The idea of the proof is to eliminate all the terms in (4.10), except the
first one,
using (5.4) so that the corresponding curvature terms are added. As a first
step
substituting
$A_{ijlkm}=A_{ijklm}+R_{ijkl,m}$ to (4.10) we obtain
$$2A_{ijklm}+R_{ijkl,m}+A_{ijmkl}+A_{ikljm}+A_{ikmjl}+A_{ilmjk}=0\;.\leqno
(5.21)$$
Now using the symmetry relations in Proposition 4.2(i), we can write
$$A_{ijmkl}=R_{ijkm,l}+A_{ijkml}=R_{ijkm,l}+A_{ijklm}\;,\leqno (5.22)$$
so (5.21) becomes
$$3A_{ijklm}+R_{ijkl,m}+R_{ijkm,l}+A_{ikljm}+A_{ikmjl}+A_{ilmjk}\;.\leqno
(5.23)$$
In the same way we can treat all components of the normal tensors in (5.21)
such that
either $j$ or $k$ (but not both) occurs as the  the second or third subscript.
Eliminating in this way $A_{ikljm}$ and $A_{ikmjl}$ from (5.23) we obtain
$$5A_{ijklm}+R_{ijkl,m}+R_{ijkm,l}+R_{ikjl,m}+R_{ikjm,l}+A_{ilmjk}=0\;.\leqno (5
.24)$$
To eliminate he last term we should do the same operation twice to put both
$j$ and $k$
to the second and third place:
$$A_{ilmjk}=R_{iljm,k}+A_{iljmk}=R_{iljm,k}+A_{ijlkm}=R_{iljm,k}+R_{ijkl,m}+A_{i
jklm}\;.$$
Substituting this expression for $A_{ilmjk}$ into (5.24) we arrive to
(5.20). $\ecarre$

\ms
{\bf Proposition 5.15.} (The second Bianchi identity) {\it On any Fedosov
manifold}
$$R_{ij(kl,m)}:=R_{ijkl,m}+R_{ijlm,k}+R_{ijmk,l}=0\;.\leqno (5.25)$$

\ms
{\bf Proof.} This identity is again true in a much more general context. In
particular,
we can deduce it from the similar identity for arbitrary symmetric connections
(with the subscript $i$ lifted -- see e.g. (49.13) in [Th], p.132). But for
the convenience of
the reader we will reproduce the proof from [Th]. Another proof can be
obtained by use of normal
coordinates as in [ON].

Let us interchange $l$ and $m$ in (5.20) and subtract the result from
(5.20) using
the identities obtained from the identities (a), (b), (c) from Theorem 5.13
by taking
the first covariant derivative. We will get then:
$$0=2R_{ijkl,m}-2R_{ijkm,l}+R_{ijkm,l}-R_{ijkl,m}
+R_{ikjl,m}-R_{ikjm,l}+$$
$$+R_{ikjm,l}-R_{ikjl,m}+R_{iljm,k}-R_{imjl,k}=$$
$$R_{ijkl,m}-R_{ijkml}+R_{iljm,k}-R_{imjl,k}=$$
$$R_{ijkl,m}+R_{ijmk,l}+(-R_{ilmj,k}-R_{imjl,k})=$$
$$R_{ijkl,m}+R_{ijmk,l}+R_{ijlm,k}\;.\qquad \ecarre$$

\ms
{\bf Proposition 5.16.} (Integrability identity.) {\it On any Fedosov manifold}
$$R_{imkj,l}+R_{ijml,k}+R_{iljk,m}+R_{iklm,j}=0\;.\leqno (5.26)$$
\ms
{\bf Proof.} Let us use the identity (4.18) with $r=2$ i.e.
$$A_{ijklm}-A_{kjilm}-A_{ilkjm}+A_{klijm}=0\;.\leqno (5.27)$$
Substituting here the expressions  of $A_{ijklm}$ and other terms
through the components $R_{ijkl,m}$ given by Lemma 5.14, and using (2.4),
(2.7) we obtain:
$$-6(A_{ijklm}-A_{kjilm}-A_{ilkjm}+A_{klijm})=2(R_{ijkl,m}+R_{jkli,m}+R_{klij,m}
+R_{lijk,m})+$$
$$+(R_{jikm,l}+R_{jkmi,l})+(R_{iljm,k}+R_{ijml,k})+
(R_{limk,j}+R_{lkim,j})+(R_{klmj,i}+R_{kjlm,i})\;.$$
The expression in the first parentheses in the right hand side vanishes
because
the differentiated identity (2.9) gives $R_{(ijkl),m}=0$. Now using the
differentiated
first Bianchi identity (2.6) we can rewrite the right hand side as
$$-(R_{milj,k}+R_{mjik,l}+R_{mkjl,i}+R_{mlki,j})\;.$$
This expression should vanish due to (5.27). Interchanging $m$ and $i$,
as well as $k$ and $l$, we arrive to (5.26). $\ecarre$

\ms
{\bf Remark 5.17.} The combinatorial structure of the integrability
identity (5.25)
has the following interesting features.

\ms
1) The first subscript $i$ is fixed, but each of other subscripts $j,k,l,m$
occurs exactly once
in any of the four available positions (second, third, fourth and fifth) in
$R_{ijkl,m}$.

\ms
2) In the cyclic order of positions ($2\to 3\to 4\to 5\to 2$) each of the
subscripts
$j,k,l,m$ is followed exactly twice by another subscript and exactly once
by two others.
In this way $j$ is followed twice by $k$, $k$ by $j$, $l$ by $m$ and $m$ by
$l$,
so the preference of double following is mutual and splits the subscripts
$j,k,l,m$
in two pairs: $\{j,k\}$ and $\{l,m\}$.

\bs
Now we can prove that already established identities for $R_{ijkl,m}$ form
a complete set
at a point.

\ms
{\bf Theorem 5.18.} {\it Assume that we are given a tensor $R_{ijkl,m}$ at
the origin
$0\in\R^n$. Then  necessary and sufficient conditions for this tensor to be
the first
covariant derivative of the curvature
tensor of  a Fedosov manifold at a point are:}
$$ R_{ijkl,m}=-R_{ijlk,m}\;; \leqno (\hbox{a}_1)$$
$$ R_{ijkl,m}=R_{jikl,m}\;; \leqno (\hbox{b}_1)$$
$$R_{i(jkl),m}:=R_{ijkl,m}+R_{iklj,m}+R_{iljk,m}=0\;;\leqno (\hbox{c}_1)$$
$$R_{ij(kl,m)}:=R_{ijkl,m}+R_{ijlm,k}+R_{ijmk,l}=0\;.\leqno (\hbox{d}_1)$$
$$R_{imkj,l}+R_{ijml,k}+R_{iljk,m}+R_{iklm,j}=0\;.\leqno (\hbox{e}_1)$$

\ms
{\bf Proof.} The conditions (a$_1$)--(e$_1$) are necessary
as was established above.
Let us prove that they are sufficient. This is done similarly to the proof of
Theorem 5.13.

Let us introduce a  hypothetical normal
tensor $A_{ijklm}$ at $0$ by the formula (5.20). We should first check
that the tensor $A_{ijkl}$ satisfies the conditions
from the part 2) of Theorem 4.5 with $r=2$ i.e. the symmetries
$$A_{ijklm}=A_{ikjlm}, \quad A_{ijklm}=A_{ijkml} \;,\leqno (5.28)$$
the identity (4.10) and the identity (5.27).
This can be done by substitution of the expression (5.20) and similar one
with permuted subscripts into the desired identities and subsequent use
of the properties (a$_1$)--(e$_1$).  For example, to prove the first
equality in (5.28) interchange $j$ and $k$ in (5.20) and subtract the right
hand sides. We get then
$$-6(A_{ijklm}-A_{ikjlm})=R_{ijkl,m}-R_{ikjl,m}+R_{iljm,k}-R_{ilkm,j}=$$
$$(R_{ijkl,m}+R_{iklj,m})-R_{ilmj,k}-R_{ilkm,j}=-R_{iljk,m}-R_{ilmj,k}-R_{ilkm,j
}=0\;,$$
where we used (a$_1$), (c$_1$) and (d$_1$). Similar straightforward
calculations
lead to the second equality in (5.28),
as well as to (4.10). The proof of (5.27) can be obtained by reversing the
calculations in
the proof of Proposition 5.16 if we note first that the ``differentiated"
identity (2.9)
follows from (a$_1$)--(c$_1$) as in the proof of Proposition 2.4.

Now we can supplement the obtained tensors $A_{ijklm}$ by arbitrary set of
tensors
$A_{ijk\al_1\dots\al_r},\ r\ne 2,$ satisfying the conditions of the second
part of
Theorem 4.5 (for example we can take $A_{ijk\al_1\dots\al_r}=0,\ r\ne 2$).
By Theorem 4.5
there exists a structure of Fedosov manifold near $0$ in $\R^n$ with the
normal tensors
prescribed in this way.

It remains to check that the obtained Fedosov manifold indeed has the
prescribed
first covariant derivative of the curvature tensor at $0$ i.e. that (5.4) holds
with $A_{ijklm}$ as constructed above and with
$R_{ijkl,m}$ coinciding with the given initial tensor. This is also
straightforward
(the properties (a$_1$)--(d$_1$) should be used).
$\ecarre$

\ms
{\bf Remark 5.19.} Clearly we can prescribe arbitrarily both tensors
$R_{ijkl}$ and
$R_{ijkl,m}$ at $0$ provided they satisfy the conditions of Theorems 5.13
and 5.18.

\vfill\eject
\centerline{\bf References.}


\bs
\item{[BFFLS]} F.Bayen, M.Flato, C.Fronsdal, A.Lichnerowicz, D.Sternheimer:
Deformation theory and quantization.I. Ann. Phys., {\bf 111} (1978), 61-110

\ms
\item{[Br]} R.L.Bryant: An introduction to Lie groups and symplectic
geometry. In: Geometry and Quantum Field Theory. D.S.Freed, K.Uhlenbeck eds.,
IAS/Park City Mathematics Series, v.1. Amer. Math. Society, Institute for
Advanced Study,
1995, 7-181

\ms
\item{[C1]} E.Cartan: Sur une g\'eneralisation de la notion de
courbure de Riemann et les espaces \'a torsion. C.R. Acad. Sci.
Ser. A {\bf 174} (1922), 593-597

\ms
\item{[C2]} E.Cartan: Sur les vari\'etes \'a connexion affine
et la th\'eorie de la relativit\'e g\'en\'eralis\'ee. Ann. Ec. Norm.
Sup. {\bf 40} (1923), 325-412

\ms
\item{[C3]} E.Cartan: On manifolds with an affine connection
and the theory of general relativity. English translation with
a foreword by A.Trautman. Bibliopolis, Napoli, 1986

\ms
\item{[D]} P.Deligne: D\'eformations de l'Algebre des Fonctions d'une
Vari\'ete Symplectique: Comparison entre Fedosov et De Wilde, Lecompte.
Selecta Mathematica, New Series, {\bf 1} (1995), 667-697

\ms
\item{[DW-L]} M.De Wilde. P.B.A.Lecompte: Existence of star-products
on exact symplectic manifolds. Annales de l'Institut Fourier
{\bf 35} (1985), 117-143

\ms
\item{[D-F-N]} B.A.Dubrovin, A.T.Fomenko, S.P.Novikov:
Modern Geometry - Methods and Applications. Part I.
Springer-Verlag, 1984, 1992

\ms
\item{[E]} L.P.Eisenhart: Non-riemannian geometry. American Mathematical
Society, 1927

\ms
\item{[F1]} B.V.Fedosov: A simple geometrical construction of deformation
quantization.
J. Diff.Geom., {\bf 40} (1994), 213-238

\ms
\item{[F2]} B.V.Fedosov: Deformation quantization and index theory.
Akademie Verlag, Berlin, 1996

\ms
\item{[Ha]} K.Habermann: Basic properties of symplectic Dirac operators.
Commun. Math. Phys., {\bf 184}(1997), 629-652
\ms
\item{[He]} S.Helgason: Differential Geometry, Lie Groups, and Symmetric
Spaces.
\smallskip
Academic Press, 1978

\ms
\item{[Hs]} H.Hess: Connections on symplectic manifolds and geometric
quantization.
In: Differential Geometric Methods in Math. Physics.
Proc. Conf. Aix-en-Province and Salamanca, 1979.
Lect. Notes in Math., {\bf 836}, 153-166.
Springer-Verlag, 1980

\medskip
\item{[K-N]} S.Kobayashi, K.Nomizu: Foundations of Differential Geometry.
\smallskip
Interscience Publishers, New York, Volume I, 1963; Volume II, 1969

\ms
\item{[Ko]} M.Kontsevich: Deformation quantization. Preprint, 1997


\ms
\item{[Le1]} H.-C.Lee: A kind of even-dimensional geometry and its application
to exterior calculus. Amer. Journal Math., {\bf 65} (1943), 433-438

\ms
\item{[Le2]} H.-C.Lee: On  even-dimensional skew-metric spaces and their
group of
transformations. Amer. Journal Math., {\bf 67} (1945), 321-328

\ms
\item{[Lm]} V.G.Lemlein: On spaces with symmetric almost symplectic connection.
Doklady Akademii nauk SSSR, {\bf 115}, no.4 (1957), 655-658 (in Russian)

\ms
\item{[Lv]} T.Levi-Civita: Nozione di parallelismo in una varieta
qualunque e consequente specificazione geometrica della curvature
Riemanniana. Rend. Palermo, {\bf 42} (1917), 173-205

\ms
\item{[Li1]} P.Libermann: Sur le probl\`eme d'\'equivalence de certaines
structures
infinit\'esimales. Ann. Mat. Pura Appl., {\bf 36} (1954), 27-120

\ms
\item{[Li2]} P.Libermann: Sur les structures presque complexes et autres
structures
infinit\'esimales r\'eguli\`eres. Bull. Soc. Math. France, {\bf 83} (1955),
195-224

\ms
\item{[Li3]} P.Libermann: Probl\`emes d'\'equivalence et g\'eom\'etrie
symplectique.
Ast\'erisque, {\bf 107-108} (1983), 43-68

\ms
\item{[Lch]} A.Lichnerowicz: Quantum mechanics and deformations of
geometric dynamics.
In: Quantum theory, Groups, Fields and Particles. (Barut, A.O., ed.), 3-82.
Reidel, Dordrecht, 1983

\ms
\item{[M]} J.M.Morvan: Quelques invariants topologiques en g\'eom\'etrie
symplectique.
Ann. Inst. H.Poincar\'e, {\bf 38} (1983), 349-370

\ms
\item{[ON]} B.O'Neill: Semi-Riemannian geometry with applications to
relativity.
Academic Press, 1983

\ms
\item{[S]} R.W.Sharpe: Cartan's generalization of Klein's
Erlangen Program. Springer, 1997

\ms
\item{[Ta]} D.E.Tamarkin: Topological invariants of connections on
symplectic manifolds.
Functional Analysis and its Applications, {\bf 29}, no.4 (1995), 45-56

\ms
\item{[Th]} T.Y.Thomas: The differential invariants of generalized spaces.
Cambridge University Press, London, 1934

\ms
\item{[To]} Ph.Tondeur: Affine Zusammenh\"ange auf Mannigfaltigkeiten mit
fast-symplectischer Struktur. Comment. Math. Helvetici, {\bf 36} (1961),
234-244

\ms
\item{[Va1]} I.Vaisman: Symplectic curvature tensors. Monatshefte f\"ur
Mathematik,
{\bf 100}, no.4 (1985), 299-327

\ms
\item{[Va2]} I.Vaisman: Lectures on the geometry of Poisson manifolds.
Birkh\"auser, 1994

\ms
\item{[V]} O.Veblen: Invariants of quadratic differential forms.
Cambridge tracts in Math. and Math. Physics, 24. Cambridge Univ. Press, 1927;
reprinted 1962

\ms
\item{[V-T]} O.Veblen, T.Y.Thomas: The geometry of paths. Transactions of
Amer. Math. Soc.,
{\bf 25} (1923), 551-608

\ms
\item{[Ve]} J.Vey: D\'eformation du crochet de Poisson sur une vari\'et\'e
symplectique.
Comment. Math. Helv., {\bf 50} (1975), 421-454

\ms
\item{[W]} A.Weinstein: Deformation quantization. S\'eminaire Bourbaki,
juin, Paris 1994

\ms
\item{[We]} H.Weyl: Reine Infinitesimalgeometrie. Math. Z. {\bf 2} (1918),
384-411
\end